\documentclass[reqno]{amsart}
\usepackage{amsmath}
\usepackage{amscd}
\usepackage{amsfonts}
\usepackage{enumerate}
\usepackage{comment}

\usepackage[plainpages=false,colorlinks,hyperindex,pdfpagemode=None,bookmarksopen,linkcolor=red,citecolor=blue,urlcolor=blue]{hyperref}
\usepackage{pdflscape}

\usepackage[a4paper,left=12mm,right=12mm,top=15mm,bottom=15mm]{geometry}

\newtheorem{thm}{Theorem}%[section]
\newtheorem{lem}[thm]{Lemma}
\newtheorem{eg}[thm]{Example}
\newtheorem{prop}[thm]{Proposition}

\newtheorem{defn}[thm]{Definition}

\newcommand{\BC}{\mathbb{C}}
\newcommand{\BN}{\mathbb{N}}
\newcommand{\BR}{\mathbb{R}}
\newcommand{\RP}{\mathbb{R}_+}

\newcommand{\CL}{\mathcal{L}}

\newcommand{\CP}{\mathrm{CP}}

\newcommand{\Tr}{\mathrm{Tr}}
\newcommand{\tp}{\mathrm{t}}

\newcommand{\R}{\mathrm{Re}\ \!}

\newcommand{\sa}{\mathrm{sa}}
\newcommand{\id}{\mathrm{id}}

\newcommand{\smnoind}{\smallskip\noindent}

\begin{document}
\title{Non-unital operator systems that are dual spaces}

\author{Yu-Shu Jia and Chi-Keung Ng}

\address[Yu-Shu Jia]{Chern Institute of Mathematics, Nankai University, Tianjin 300071, China.}
\email{1410035@mail.nankai.edu.cn}

\address[Chi-Keung Ng]{Chern Institute of Mathematics and LPMC, Nankai University, Tianjin 300071, China.}
\email{ckng@nankai.edu.cn}

\date{\today}


\begin{abstract}
We will extend the main result of \cite{dual} to the non-unital case with a totally different proof; i.e., we give an abstract characterization of an arbitrary self-adjoint weak$^*$-closed subspace of $\CL(H)$ (equipped with the induced matrix norm, the induced matrix cone and the induced weak$^*$-topology).
In order to do this, we obtain a matrix analogues of a result of Bonsall for $^*$-operator spaces equipped with closed matrix cones. 
On our way, we observe that for a $^*$-vector $X$ equipped with a matrix cone (in particular, when $X$ is an operator system or the dual space of an operator system), a linear map $\phi:X\to M_n$ is completely positive if and only if linear functional $[x_{i,j}]_{i,j}\mapsto \sum_{i,j=1}^n \phi(x_{i,j})_{i,j}$ on $M_n(X)$ is positive. 
\end{abstract}

\maketitle

\section{Introduction}

\medskip

The aim of the article is to give the non-unital version of the main result of \cite{dual} as follows: 

\medskip

\begin{thm}[Blecher-Magajna]\label{thm:orig-BM}
Let $S$ be a unital operator system (see Definition \ref{def:os}(a)) such that it is also a dual operator space. 
Then $S$ is a dual unital operator system (see Definition \ref{def:os}(b)).  
\end{thm}

\medskip

Notice that in the unital case, the involution on $S$ is automatically weak$^*$-continuous (i.e., $S_\sa$ is weak$^*$-closed) and the matrix cone $(M_n(S)_+)_{n\in \BN}$ of $S$ is automatically weak$^*$-closed.
In fact, these two facts hold because both $(M_n(S)_+)_{n\in \BN}$ and the closed unit ball of $S_\sa$ are determined by the matrix norm and the order unit of $S$ (together with the weak$^*$-closedness of the norm-closed unital ball of $S$ as well as the Krein-Smulian theorem). 
However, when $S$ is a non-unital operator system and is a dual operator space, there is no guarantee that any of the above two facts hold. 
Indeed, we found an example where the matrix cone is not weak$^*$-closed (Example \ref{eg:counter}).
For this reason, in our generalization of Theorem \ref{thm:orig-BM}, we need to add the assumptions that the involution is weak$^*$-continuous and the matrix cone is weak$^*$-closed (Theorem \ref{thm:ext-BM}). 

\medskip

As noted in the above, some techniques in the study of unital operator systems do not work in the non-unital case. 
Therefore, we will give a totally different proof for Theorem \ref{thm:ext-BM}. 
The main tools are some results from \cite{Ng2} and \cite{Ng1}, as well as the following fact that we will verify in Proposition \ref{main}:
\begin{quotation}
Let $V$ be an operator space equipped with a completely isometric involution and a closed matrix cone, the positive part $B_{M_n(V)}^+$ of the closed unit ball of $M_n(V)$ is weak$^*$-dense in $B_{M_n(V^{**})}^+$, where $V^{**}$ is equipped with the bidual matrix cone and the weak$^*$-topology on $M_n(V^{**})$ is as defined in Definition \ref{def:weak-st-top}. 
\end{quotation}
Note that in the above, we consider $B_{M_n(V^{**})}^+$ but not $B_{M_n(V)^{**}}^+$ (otherwise, the result is well-known). 
Moreover, for our application in Theorem \ref{thm:ext-BM}, we need to consider the case when the cone on $V$ may not even be proper (i.e. $V^+ \cap - V^+$ may be non-zero), and hence $V$ cannot be an operator system.  

\medskip

In the case of  ordered Banach spaces, one way to prove the corresponding fact of the above displayed statement is to empoly the following result of Bonsall (see \cite[Theorem 1.1.1]{wong}). 
Recall that a \emph{cone} $C\subseteq E$ is a subset satisfying $C+C\subseteq C$ and $\RP  \cdot C\subseteq C$, and $q:C\to \BR$ is \emph{sublinear} if $q(x+y)\leq q(x) + q(y)$ and $q(t\cdot x) = t q(x)$ ($x,y\in C$ and $t\in \RP$). 

\medskip

\begin{prop}[Bonsall]\label{Bonsall}
	Let $E$ be a real ordered vector space with a cone $E_+$. 
	Suppose that $p$ is a sublinear function from $E$ to $\BR$. 
	If $q$ is a sublinear function from $E_+$ to $\BR$ with  	
	\[-q(u)\leq p(u)\qquad (u\in E_+).\]
	Then one can find a linear functional $g$ on $E$ such that
	\begin{equation*}
		-q(u)\leq g(u) \quad \text{and}\quad g(x)\leq p(x) \qquad (u\in E_+;x\in E).
	\end{equation*}
\end{prop}

\medskip

Therefore, we will first give an analogue of this result in the matrix case (see Lemma \ref{lem:mat-Bonsall}), and use it to obtain the above displayed statement. 
On our way to this matrix analogue, we observe that (see the proof of Lemma \ref{bijective}): 
\begin{quotation}
Let $n\in \BN$ and $X$ be a semi-matrix ordered vector space, in the sense of Definition \ref{MOS}(a) (in particular, an operator system or the dual space of an operator system). 
The assignment $\phi\mapsto \Upsilon_{\phi}$ is a linear bijection from the space of linear maps from $X$ to $M_n$, to the space of linear functionals on $M_n(X)$ such that $\phi$ is completely positive if and only if $\Upsilon_{\phi}$ is positive; where
$$\Upsilon_{\phi}\big([x_{k,l}]_{k,l}\big):= {\sum}_{k,l=1}^n \phi(x_{k,l})_{k,l}\qquad \big([x_{k,l}]_{k,l} \in M_n(X)\big).$$
\end{quotation}

\medskip

\section{Notations}\label{sec:notation}

\medskip

For $m,n\in \BN$, we denote by $M_n$ and $M_{n,m}$ the spaces $M_n(\mathbb{C})$ and $M_{n,m}(\mathbb{C})$ of complex matrices. 
Suppose that $e_k:=(0,\cdots,0,1,0,\cdots,0)\in M_{1,n}$ is 
the matrix with only the $k^\mathrm{th}$-entry being $1$.
We set 
\begin{equation}\label{eqt:def-tau-kl}
	e^{k,l}:= e_k^*e_l\in M_n \quad \text{and} \quad \beta_e:= (e_1, \dots, e_n)\in M_{1,n}(M_{1,n}).
\end{equation}
For a complex vector space $X$, we denote $M_n(X) :=  M_n\otimes X$, and see elements in $M_n(X)$ as matrices with entries being elements of $X$. 
In this case, 
\begin{equation}\label{eqt:id-tensor}
\alpha\otimes x \quad \text{is identified with}\quad  [\alpha_{i,j}x]_{i,j}, \quad \text{for } \alpha = [\alpha_{i,j}]_{i,j}\in M_n \text{ and } x\in X.
\end{equation}
Under this identification, one has
\begin{equation}\label{eqt:beta-alpha}
	\beta_e(\alpha\otimes \gamma) \beta_e^* = {\sum}_{i,j=1}^n \alpha_{i,j}\gamma_{i,j}\qquad  (\alpha, \gamma\in M_n). 
\end{equation}
Moreover, we will regard $M_n(X)\subseteq M_{n+1}(X)$ by putting elements of $M_n(X)$ into the upper left corner of $M_{n+1}(X)$. 
%Under this identification, we set
%\[M_\infty(X):={\bigcup}_{n\in \mathbb{N}}M_n(X).\]
For $x\in M_n(X)$ and $y\in M_m(X)$, we write
$x\oplus y:=\begin{pmatrix}
	x&0\\
	0&y
\end{pmatrix}\in M_{n+m}(X).$ 
If $Y$ is another complex vector space, we denote by $L(X,Y)$ the set of complex linear maps from $X$ to $Y$. 
For $\varphi\in L(X,Y)$ and $n\in \BN$, we define $\varphi^{(n)}:M_n(X)\to M_n(Y)$ by
\[\varphi^{(n)}([x_{i,j}]_{i,j}):=[\varphi(x_{i,j})]_{i,j}\quad ([x_{i,j}]\in M_n(X)).\]
%and put $\varphi^{(n)}:=\varphi^{(\infty)}|_{M_n(X)}$.

\medskip

We will need the notions of operator spaces, dual operator spaces, completely bounded maps, complete contractions and  complete isometries. 
Reader who is not familiar with these concepts may consult \cite{OPS}. 

\medskip

\begin{defn}\label{MOS}
	(a) A complex vector space $X$ is called a \emph{semi-matrix ordered vector space} if there is an involution $^*:X\to X$ and for each $n\in \BN$, there is a subset $M_n(X)_+$  of $M_n(X)_{sa}:=\left\lbrace x\in M_n(X):x^*=x\right\rbrace$
such that for $m,n\in \BN$, 
\begin{itemize}
	\item $\gamma^*v\gamma\in M_n(X)_+$ whenever $v\in M_m(X)_+$ and $\gamma\in M_{m,n}$;
	\item $u\oplus v\in M_{m+n}(X)_+$ whenever $u\in M_m(X)_+$ and $v\in M_n(X)_+$. 
\end{itemize}	In this case, $(M_n(X)_+)_{n\in \BN}$ is called the \emph{matrix cone} of $X$.

\smnoind
(b) Let $X$ and $Y$ be semi-matrix ordered vector spaces and $\varphi\in L(X,Y)$. 
We define $\varphi^*\in L(X,Y)$ by
\begin{equation}\label{eqt:def-map-invol}
\varphi^*(x):= \varphi(x^*)^*\qquad (x\in X). 
\end{equation}
Then $\varphi$ is said to be \emph{completely positive} if $\varphi^*=\varphi$ and $\varphi^{(n)}(M_n(X)_+)\subseteq M_n(Y)_+$ ($n\in \BN$). 
We will denote by $\CP(X,Y)$ the set of all completely positive maps and set 
$L(X,Y)_\sa:= \{\varphi\in L(X,Y): \varphi^* = \varphi\}$.

\smnoind
(c) If $X$ is both an operator space and a semi-matrix ordered vector space such that the induced involution $^*:M_n(X)\to M_n(X)$ is isometric and $M_n(X)_+$ is norm-closed  for each $n\in \BN$, then $X$ is called a \emph{semi-matrix ordered operator space (SMOS)}. 
\end{defn}

\medskip

Notice that unlike \cite{Werner}, we DO NOT assume $M_n(X)_+\cap -M_n(X)_+ = \{0\}$ and hence, we added the prefix ``semi-'' in our definitions.

\medskip

\section{Main Result}\label{section2}	

\medskip
We denote by $\Tr_n: M_n \to \BC$ the usual trace on $M_n$; i.e. $\Tr_n\big([\beta_{i,j}]_{i,j}\big) := \sum_{i=1}^n \beta_{i,i}$. 
For $\tau\in M_m(M_n)$, we consider $\tau^\tp \in M_m(M_n)$ to be the element satisfying 
$$(\tau^\tp)_{k,l} = (\tau_{l,k})^\tp \in M_n\qquad (k,l=1,\dots,m),$$ 
where $(\tau_{l,k})^\tp$ is the usual transpose in  $M_n$.

\medskip

If $m,n\in \BN$, we may identify $M_m(M_n) = M_m\otimes M_n$ with $M_{m\times n}$ as $C^*$-algebras in the canonical way. 
It is easy to see that $\tau\mapsto \tau^\tp$ is an order isomorphism from $M_m(M_n)$ onto itself. 
Moreover, for any $\tau \in M_m(M_n)$, the value $\Tr_m\big(\Tr_n^{(m)}(\tau)\big)$ equals the trace of $\tau$, when $\tau$ is considered as an element in $M_{m\times n}$.
We set 
\begin{equation}\label{eqt:def-norm-1}
\|\tau\|_1 := \Tr_m\big(\Tr_n^{(m)}(|\tau|)\big) \qquad (\tau \in M_m(M_n)).
\end{equation}
The following result follows from the above observation. 

\medskip

\begin{lem}\label{lem:dual}
Let $m,n\in \BN$. 
	For every $\tau\in M_m(M_n)$, we defined $\varphi_\tau: M_m(M_n)\to \BC$ by $$\varphi_\tau(\alpha) := \Tr_m\big(\Tr_n^{(m)}(\tau^\tp\alpha)\big)
	\qquad (\alpha\in M_m(M_n)).$$
	Then $\tau\mapsto \varphi_\tau$ is an isometric order isomorphism from $\big(M_m(M_n), \|\cdot\|_1\big)$ onto the dual space of the $C^*$-algebra $M_m(M_n)$.
\end{lem}

\medskip

When $L(M_m,M_n)$ is equipped with the involution as in \eqref{eqt:def-map-invol} as well as 
the cone $\CP(M_m,M_n)$ (see Definition \ref{MOS}(b)), we know from  \cite[Lemma 2.1]{choi} that the map 
$$\theta:M_{m}(M_n)\to L(M_m,M_n)$$ 
given by 
$\theta_\tau(\alpha):={\sum}_{k,l=1}^m\tau_{k,l}\cdot \alpha_{k,l}$ ($\tau = [\tau_{i,j}]_{i,j} \in M_{m}(M_n); \alpha=[\alpha_{i,j}]_{i,j}\in M_m$)
is an order isomorphism. 
Observe that
\begin{equation}\label{eqt:theta-tau-otimes}
	(\theta_\tau\otimes \id_X)(v) = {\sum}_{k,l=1}^m\tau_{k,l}\otimes v_{k,l} \qquad \big(\tau= [\tau_{i,j}]_{i,j}\in M_m(M_n); v=[v_{i,j}]_{i,j}\in M_m(X)\big).
\end{equation}
Our next result can be regarded as the dual version of \cite[Lemma 2.1]{choi}.
However, since this result concern with the more general case of semi-matrix ordered vector spaces instead of operator systems, one cannot use \cite[Lemma 2.1]{choi} to obtain this result, even if it is known that $L(L(M_n,X),\BC)$ is order isomorphic to $L(X, M_n)$. 

\medskip

\begin{lem}\label{bijective}
Suppose that $X$ is a semi-matrix ordered vector space and $n\in \BN$.  
If we equip $L(X, M_n)$ with the involution as in \eqref{eqt:def-map-invol} and the cone $\CP(X,M_n)$, then $\Theta: F\mapsto \Theta_F$ is an order isomorphism from $L(M_n(X),\BC)$ onto $L(X, M_n)$; where the linear map $\Theta_F:X\to M_n$ is defined by 
	\begin{equation}\label{eqt:def-Delta}
		\Theta_F(x)_{k,l}:=F(e^{k,l}\otimes x) \qquad (x\in X; k,l=1,\dots,n).
	\end{equation}
\end{lem}
\begin{proof}
Clearly, $\Theta$ is a complex linear injection. 
	Moreover, for any $F\in L(M_n(X),\BC)_\sa$, it follows from 
	$$\Theta_F(x^*)_{k,l} = F(e^{k,l}\otimes x^*) = F((e^{l,k}\otimes x)^*) = F(e^{l,k}\otimes x)^* = (\Theta_F(x)_{l,k})^*\qquad (x\in X; k,l\in \{1,\dots,n\})$$
	that $\Theta_F\in L(X,M_n)_\sa$. 
	In other words, $\Theta$ preserves the involutions. 
	
	Let us consider a map $\Upsilon: \phi\mapsto \Upsilon_\phi$ from $L(X,M_n)$ to $L(M_n(X),\BC)$ given by 
	\begin{equation*}
		\Upsilon_\phi(u):={\sum}_{k,l=1}^n\phi(u_{k,l})_{k,l} \qquad (u\in M_n(X)). 
	\end{equation*}
Pick any $\phi\in L(X, M_n)$. 
One has  $\Theta_{\Upsilon_{\phi}}(x)_{k,l} = \Upsilon_\phi(e^{k,l}\otimes x) = \phi(x)_{k,l}$. 
Hence, $\Theta$ is surjective and $\Upsilon$ is the inverse of $\Theta$. 

Moreover, by Relation \eqref{eqt:beta-alpha}, we know from $\phi^{(n)}(u) = \sum_{k,l=1}^n e^{k,l}\otimes \phi(u_{k,l})$ that 
\begin{equation}\label{eqt:def-Upsilon-phi}
\Upsilon_\phi(u) = {\sum}_{k,l=1}^n\phi(u_{k,l})_{k,l}=\beta_e\phi^{(n)}(u)\beta_e^*\qquad (u\in M_n(X)).
\end{equation}
From this, we see that $\Upsilon$ is positive.

Now, consider $F\in L(M_n(X),\BC)$, $m\in \BN$, $w\in M_m(X)$ and $\tau\in M_m(M_n)$. 
Then Relations \eqref{eqt:beta-alpha},  \eqref{eqt:theta-tau-otimes} and \eqref{eqt:def-Upsilon-phi} give
\begin{align}\label{eqt:Upsil-phi}
	F\big((\theta_\tau\otimes \id_X)(w)\big) 
	&= \Upsilon_{\Theta_F}\big((\theta_\tau\otimes \id_X)(w)\big) 
	= \beta_e \Theta_F^{(n)}\big({\sum}_{k,l=1}^m \tau_{k,l}\otimes w_{k,l}\big) \beta_e^*
	=  \beta_e \big({\sum}_{k,l=1}^m \tau_{k,l}\otimes \Theta_F(w_{k,l})\big) \beta_e^*\nonumber\\
	&= {\sum}_{k,l=1}^m {\sum}_{i,j=1}^n (\tau_{k,l})_{i,j} \Theta_F(w_{k,l})_{i,j}
	= \Tr_m\big(\Tr_n^{(m)}\big(\tau^\tp\Theta_F^{(m)}(w)\big)\big).
\end{align}

In order to verify that $\Theta$ is positive, let us assume that $F\in L(M_n(X),\BC)_+$. 
As $\Theta$ preserves the involutions, one has $\Theta_F\in L(X,M_n)_\sa$. 
	By Lemma \ref{lem:dual}, in order to establish $\Theta_F^{(m)}(w)\geq 0$, it suffices to show that for each $\tau \in M_m(M_n)_+$, 
	$$\Tr_m\big(\Tr_n^{(m)}\big(\tau^\tp\Theta_F^{(m)}(w)\big)\big)\geq 0;$$ 
	or via Relation \eqref{eqt:Upsil-phi}, $F\big((\theta_\tau\otimes \id_X)(w)\big) \geq 0$. 
	Indeed, let $\tau \in M_m(M_n)_+$. 
	Since $\theta$ is an order isomorphism (see \cite[Lemma 2.1]{choi}), we know that $\theta_\tau$ is completely positive.
	Thus, \cite[Lemma 5.1.6]{OPS} and \cite[Theorem 1]{choi2} imply that there exist $\gamma_1,\dots, \gamma_k\in M_{m,n}$ with
	$\theta_\tau(\alpha)={\sum}_{i=1}^k \gamma_i^*\alpha \gamma_i$ ($\alpha\in M_m$). 
	Therefore, 
	$$(\theta_\tau\otimes \id_X)(w) ={\sum}_{i=1}^k \gamma_i^*w \gamma_i\in M_n(X)_+,$$
	and one has $F((\theta_\tau\otimes \id_X)(w)) \geq 0$, as required. 
\end{proof}

\medskip

\begin{defn}(\cite{effros1})
	A \emph{finite matrix gauge} on a complex vector space $X$ is a collection $(\rho_n)_{n\in \BN}$ of sublinear functional $\rho_n: M_n(X)\to \RP$ such that
	for $v\in M_m(X),w\in M_n(X)$ and $\alpha\in M_{m,n}$, one has 
	\begin{itemize}
		\item $\rho_{m+n}(v\oplus w)=\max\left\lbrace\rho_m(v),\rho_n(w)\right\rbrace$;
		\item $\rho_n(\alpha^*v\alpha)\leq\left\|\alpha\right\|^2\rho_m(v)$.
	\end{itemize}
\end{defn}

\medskip

In the following, we denote by $I_{m\times n}$ the identity of $C^*$-algebra $M_m(M_n)$. 

\medskip

We can now present the following matrix analogue of Proposition \ref{Bonsall}. 
Notice that we need to assume $\rho$ to be a finite matrix gauge instead of a ``matrix sublinear map'', and $\phi$ to be a self-adjoint linear map  from $X$ to $M_n$ instead of $-\phi$  being a ``matrix sublinear map'' from $X_+$ to $(M_n)_\sa$.

\medskip

\begin{lem}\label{lem:mat-Bonsall}
	Let $X$ be a semi-matrix ordered vector space and $\rho$ be a finite matrix gauge on $X$. 
	For $n\in \BN$, if $\phi\in  L(X,M_n)_\sa$ satisfies
	\begin{equation*}\label{eqt:phi-infty-rho}
		\phi^{(m)}(w) \leq \rho_m(w)I_{m\times n} \qquad (w\in M_m(X)_+; m\in \BN),
	\end{equation*}
	then there exists $\psi\in L(X,M_n)_\sa$ such that
	\begin{equation}\label{eqt:phi-psi}
	\phi^{(m)}(u) \leq \psi^{(m)}(u)\quad \text{and} \quad \psi^{(m)}(z) \leq \rho_m(z)I_{m\times n} \qquad  (u\in M_m(X)_+;z\in M_m(X)_\sa;m\in \BN).
	\end{equation}
\end{lem}
\begin{proof}
As in \cite[p.139]{effros1}, we define a sublinear functional $\lambda:M_n(X)\to \RP$ by 
\begin{equation*}\label{eqt:def-hat-rho}
	\lambda(y):=\inf\big\{\|\tau\|_1\rho_m(z):y=(\theta_\tau \otimes \id_X)(z), \text{ where }z\in M_m(X); \tau\in M_m(M_n)_+; m\in \BN\big\} \quad (y\in M_n(X));
\end{equation*}
here $\|\tau\|_1$ is as in \eqref{eqt:def-norm-1}. 
Set 
$$F:=\Upsilon_\phi\in L(M_n(X),\BC)$$ 
(see Relation \eqref{eqt:def-Upsilon-phi}). 
Since $\phi$ is self-adjoint, Lemma \ref{bijective} implies $F(M_n(X)_\sa)\subseteq \BR$. 
We claim that $F(u)\leq \lambda(u)$ ($u\in M_n(X)_+$). 

Suppose on the contrary that $F(v)> \lambda(v)$ for some $v\in M_n(X)_+$. 
By \cite[Lemma 6.8]{effros1}, we have 
$$\lambda(v)=\inf\big\{\Tr_n(\beta^2)\rho_n(w): \beta\in (M_n)_+\text{ is invertible and } w\in M_n(X)  \text{ such that }v=\beta w\beta\big\}.$$
Hence, we can find $w\in M_n(X)$ and an invertible matrix $\beta\in (M_n)_+$ with
$$v=\beta w \beta \quad \text{and} \quad \Tr_n(\beta^2)\rho_n(w)<F(v).$$ 
Let us write $\beta$ as $(\beta_1,\dots,\beta_n)$ with $\beta_i\in M_{n, 1}$, and set
$\tilde{\beta}:=(\beta_1^*, \dots, \beta_n^*)\in M_{1,n}(M_{1,n})$. 
If we define $\check\beta:=\tilde{\beta}^*\tilde{\beta}\in M_n(M_n)_+$, then $\|\check\beta\|_1=\Tr_n(\beta^2)$ and \eqref{eqt:id-tensor} and \eqref{eqt:theta-tau-otimes} give
\[v
=\beta w \beta
= \beta w \beta^*
=\Big[{\sum}_{k,l=1}^n\beta_{i,k}w_{k,l}\bar \beta_{j,l}\Big]_{i,j}
={\sum}_{k,l=1}^n \check\beta_{k,l}\otimes w_{k,l} 
= (\theta_{\check\beta}\otimes \id_X)(w).\]
Thus,  
\begin{equation}\label{eqt:norm-check-beta}
\|\check\beta\|_1\rho_n(w)=\Tr_n(\beta^2)\rho_n(w)< F(v) = F\big((\theta_{\check\beta}\otimes \id_X)(w)\big). 
\end{equation}
On the other hand, since $\beta^{-1}\in (M_n)_+$, we know that $w=\beta^{-1}v\beta^{-1}\in M_n(X)_+$.
Therefore, the assumption on $\phi$ and Lemma \ref{lem:dual}
imply that 
\[\Tr_m\big(\Tr_n^{(m)}\big(\check \beta^\tp\phi^{(m)}(w)\big)\big) \leq\rho_m(w)\Tr_m\big(\Tr_n^{(m)}(\check \beta^\tp)\big) = \|\check \beta^\tp\|_1\rho_m(w) = \|\check \beta\|_1\rho_m(w).\]
Consequently,  Relation \eqref{eqt:Upsil-phi} tells us that $F\big((\theta_{\check\beta}\otimes \id_X)(w)\big)\leq \|\check \beta\|_1\rho_m(w)$, which contradicts Relation (\ref{eqt:norm-check-beta}).
Therefore, the above claim is established. 

Now, by Proposition \ref{Bonsall}, we obtain a linear functional $G:M_n(X)_\sa\to \BR$ such that 
\begin{equation*}
	F(u)\leq G(u) \quad \text{and}\quad G(x)\leq \lambda(x) \qquad (u\in M_n(X)_+;x\in M_n(X)_\sa).
\end{equation*}
Consider $G_0:M_n(X)\to \BC$ to be the complexification of $G$. 
Lemma \ref{bijective} implies that $\psi:= \Theta_{G_0}\in L(X,M_n)_\sa$ and we have 
$$(\psi - \phi)^{(m)}(u)\geq 0 \qquad (u\in M_m(X)_+; m\in \BN).$$
Finally, pick any $m\in \BN$ and $z\in M_m(X)_\sa$.
As $G_0(x)\leq \lambda(x)$ for every $x\in M_n(X)_\sa$, we know  that
\begin{equation}\label{G}
	G_0\big((\theta_\tau\otimes \id_X)(z)\big)\leq \|\tau\|_1\rho_m(z) =  \Tr_m\big(\Tr_n^{(m)}(\tau^\tp)\big)\rho_m(z) \qquad (\tau\in M_m(M_n)_+).
\end{equation}
It then follows from Relation \eqref{eqt:Upsil-phi} that 
\[\Tr_m\big(\Tr_n^{(m)}\big(\tau^\tp\psi^{(m)}(z)\big)\big) \leq \Tr_m\big(\Tr_n^{(m)}(\tau^\tp)\big)\rho_m(z) \qquad (\tau\in M_m(M_n)_+),\]
which implies the required inequality $\psi^{(m)}(z)\leq \rho_m(z) I_{m\times n}$, because of Lemma \ref{lem:dual}. 
\end{proof}

\medskip

In order to finish our proof, we also need the matrix version of the separation theorem.

\medskip

\begin{defn}(\cite{effros1})
Suppose $V$ is a SMOS. 
A collection $(K_n)_{n\in \BN}$ of non-empty convex sets $K_n\subseteq M_n(V)$ is called a \emph{matrix convex set}  in $V$, if for every $m,n\in \BN$, 
	\begin{itemize}
		\item $\alpha^*v\alpha\in K_n$ whenever $v\in K_m$ and $\alpha\in M_{m,n}$ satisfying $\alpha^*\alpha=I_n$;
		\item $u\oplus v\in K_{m+n}$ whenever $u\in K_m$ and $v\in K_n$. 
	\end{itemize}
Moreover, the matrix convex subset $(K_n)_{n\in \BN}$ is said to be \emph{self-adjoint} if $K_n\subseteq M_n(V)_\sa$ for each $n\in \BN$. 
\end{defn}

\medskip

For a SMOS $W$, we denote by $W^*$ the dual Banach space of $W$, and set
\[B_W:=\left\lbrace x\in W:\left\|x\right\|\leq 1\right\rbrace \quad \text{as well as} \quad B_W^+ := B_W\cap W^+.\]
For $n\in \BN$, we will identify $f=[f_{i,j}]_{i,j}\in M_n(W^*)$ with the map from $W$ to $M_n$ satisfying
$$f(x):= [f_{i,j}(x)]_{i,j} \qquad (x\in W).$$
In this case,  $f:W\to M_n$ is always completely bounded.  
In the same way, we may consider elements in $M_n(W)$ as a completely bounded map from $W^*$ to $M_n$. 

\medskip

The dual space $W^*$ is a SMOS under the dual operator space structure, the induced $^*$-operation and the matrix cone:
$$M_n(W^*)_+ := \big\{f\in M_n(W^*)_\sa: f \in \CP(W,M_n) \big\} \qquad (n\in \BN).$$

\medskip

Furthermore, we consider a duality between $M_n(W)$ and $M_n(W^*)$ via  
\begin{equation*}\label{eqt:def-pair-mat}
f(x):={\sum}^n_{k,l=1}f_{k,l}(x_{k,l}) \qquad (f\in M_n(W^*); x\in M_n(W)).
\end{equation*}

\medskip

\begin{defn}\label{def:weak-st-top}
Let $W$ be a SMOS. 

\smnoind
(a) The weakest topology on $M_n(W^*)$ under which the functional $f\mapsto f(x)$ is continuous for every $x\in M_n(W)$ will be denoted by $\sigma(M_n(W^*),M_n(W))$ and called the \emph{weak$^*$-topology} on $M_n(W^*)$. 

\smnoind
(b) A matrix convex set $(K_n)_{n\in \BN}$ in $W^*$ is said to be \emph{weak$^*$-closed} if $K_n$ is weak$^*$-closed in $M_n(W^*)$ for all $n\in \BN$. 
\end{defn}

\medskip

Observe that a net $\{f^{(\lambda)}\}_{\lambda\in \Lambda}$ in $M_n(W^*)$ converge to $f$ under the weak$^*$-topology if and only if the net $\{f_{k,l}^{(\lambda)}\}_{\lambda\in \Lambda}$ $\sigma(W^*,W)$-converges to $f_{k,l}$ in $W^*$, for all $k,l\in \{1,...,n\}$.

\medskip

Our next lemma follows more or less from the proof of \cite[Theorem 5.4]{effros1}, but we say a few words about how to apply the proof of \cite[Theorem 5.4]{effros1}. 

\medskip

\begin{lem}\label{polar}
	Let $V$ be a SMOS and $(K_n)_{n\in \BN}$ be a self-adjoint weak$^*$-closed matrix convex set in $V^{**}$ with $0\in K_1$. 
	For every $n\in \BN$ and $v_0\in M_n(V^{**})_{sa}\setminus  K_n$, there exists a completely bounded self-adjoint linear map $\phi:V\to M_n$ satisfying
$(\phi^{**})^{(m)}(w)\leq I_{m\times n}$ ($w\in K_m; m\in \BN$) but 
$(\phi^{**})^{(n)}(v_0)\not\leq I_{n\times n}$.
\end{lem}
\begin{proof}
As $0\in K_1$, we know that $0\in K_n$. 
Hence, the usual separation theorem produces a $\sigma(M_n(V^{**}),M_n(V^*))$-continuous complex linear functional $F:M_n(V^{**})\to \mathbb{C}$ such that $\R F(v)\leq 1 < \R F(v_0)$ for every $v\in K_n$. 
Since both $K_n$ and $v_0$ are self-adjoint, by replacing $F$ with $(F+F^*)/2$ if necessary, we may assume the function $F$ is self-adjoint, and have 
$$ F(v)\leq 1 < F(v_0) \qquad (v\in K_n).$$ 
By \cite[Lemma 5.3]{effros1}, there exists  a state $\omega$ on the $C^*$-algebra $M_n$ satisfying
\[ F(\alpha^* v\alpha)\leq \omega(\alpha^*\alpha)\qquad (v\in K_m, \alpha\in M_{m,n}, m\in \BN).\]
By replacing $F$ and $\omega$ with $(1-\epsilon)\ \! F$ and $(1-\varepsilon)\ \! \omega+\varepsilon\ \! \Tr_n/n$, respectively, for small enough $\epsilon>0$, we may assume that $\omega$ is faithful. 
Consider $(H, \pi, \xi_0)$ to be the GNS representation of $\omega$.
We set	
\[\bar{\alpha}:= (\alpha^*, 0, \dots, 0)^*\in M_n \quad (\alpha\in M_{1,n}),\]
and denote $\bar{H}:=\{\pi(\bar \alpha)\xi_0: \alpha\in  M_{1,n}\}$. 
The proof of \cite[Theorem 5.4]{effros1} then gives a $\sigma(V^{**},V^*)$-continuous linear map $\psi:V^{**}\to \CL(\bar{H})$ satisfying
$F(\beta^* z\alpha)=\left\langle\psi(z)\pi(\bar{\alpha})\xi_0, \pi(\bar{\beta})\xi_0\right\rangle$  ($z\in V^{**}; \alpha,\beta\in M_{1,n}$) and 
\begin{equation}\label{eqt:rel-psi}
\R \psi^{(n)}(v_0)\not\leq I_{n\times n} \quad \text{as well as}\quad \R \psi^{(m)}(w)\leq I_{m\times n}\ \ (w\in K_m; m\in \BN),
\end{equation}
when we regard $\CL(\bar{H}) = M_{\dim \bar H} \subseteq M_n$, by fixing a basis for $\bar H$. 
Observe that
\begin{align*}
\langle\psi(u)\pi(\bar{\alpha})\xi_0,\pi(\bar{\alpha})\xi_0\rangle
= F(\alpha^* u \alpha)\in \BR \qquad (u\in V^{**}_\sa; \alpha\in M_{1,n}).
\end{align*}
This implies that $\psi$ is self-adjoint and so is $\psi^{(m)}$ for all $m\in \BN$ (thus, we may remove $\R$ from \eqref{eqt:rel-psi}). 
Finally, as $\psi$ is a weak$^*$-continuous self-adjoint linear map from $V^{**}$ to $M_n$, one can find a bounded self-adjoint linear map $\phi:V\to M_n$ satisfying $\psi = \phi^{**}$.  
Since $M_n$ is finite dimensional, we know that $\phi$ is completely bounded. 
This completes the proof. 
\end{proof}

\medskip

The above actually holds when $V$ is only a $^*$-operator space, since the matrix cone plays no part in the proof.

\medskip

Now, we can proof the displayed statement in the Introduction.

\medskip

\begin{prop}\label{main}
	Suppose $V$ is a SMOS.
	Then $B_{M_{n}(V)}^+$ is weak$^*$-dense in $B_{M_{n}(V^{**})}^+$ for every $n\in \mathbb{N}$.
\end{prop}
\begin{proof}
	For $m\in \BN$, we set $K_m$ to be the weak$^*$-closure of $B_{M_m(V)}^+$ in  $M_m(V^{**})_\sa$. 
	Then $(K_m)_{m\in \BN}$ is a self-adjoint matrix convex set in $V^{**}$ with $0\in K_1$. 
	Assume on the contrary that there exists $v_0\in B_{M_n(V^{**})}^+\setminus K_n$ for some $n\in \BN$.
	By Lemma \ref{polar}, one can find a completely bounded self-adjoint linear map $\phi:V\to M_n$ such that 
$$(\phi^{**})^{(n)}(v_0)\not\leq I_{n\times n} \quad \text{and}\quad \phi^{(m)}(w)\leq I_{m\times n} \quad (w\in B_{M_m(V)}^+; m\in \BN).$$
Let us now define a finite matrix gauge $(\rho_m)_{m\in \BN}$ on $V$ by $$\rho_m(u)=\left\|u\right\| \qquad (u\in M_m(V);m\in \BN).$$ 	
The above tells us that $\phi^{(m)}(v)\leq \rho_m(v)I_{m\times n}$ for every $m\in \BN$ and $v\in M_m(V)^+$. 
By Lemma \ref{lem:mat-Bonsall}, one can find $\psi\in L(X,M_n)_\sa$  satisfying  \eqref{eqt:phi-psi}. 
Consider $m\in \BN$. 
For $w\in M_m(V)_\sa$, one knows from the second inequality of \eqref{eqt:phi-psi} that 
	$$-\|w\|I_{m\times n}\leq \psi^{(m)}(w) \leq \|w\|I_{m\times n};$$
in other words, $\|\psi^{(m)}(w)\| \leq \|w\|$. 
For $z\in M_m(X)$, 
by considering $w=\begin{pmatrix}
	0&z\\
	z^*&0
\end{pmatrix}\in M_{2m}(V)_\sa$, 
we know that 
$$\|\psi^{(m)}(z)\|\leq \max\{\|z\|, \|z^*\|\} = \|z\|. $$ 
This shows that $\psi$ is a complete contraction. 

On the other hand, since $\phi:V\to  M_n$ is completely bounded, the first inequality of \eqref{eqt:phi-psi} tells us that $\psi-\phi$ is a completely bounded completely positive map, and so is $\psi^{**} - \phi^{**}$. 
From this, we obtain the contradiction that 
$$(\phi^{**})^{(n)}(v_0)\leq (\psi^{**})^{(n)}(v_0) \leq I_{n\times n},$$
because $\psi^{**}$ is a self-adjoint complete contraction. 
\end{proof}

\medskip

\begin{defn}\label{def:os}
(a) An \emph{operator system} $S$ is a subspace of of some $\CL(H)$, equipped with the induced SMOS structure. 
Moreover, $S$ is said to be unital if it contains the identity of $\CL(H)$

\smnoind
(b) A \emph{dual operator system} $S$ is a weak$^*$-closed subspace of some $\CL(H)$, equipped with the induced SMOS structure as well as the induced weak$^*$-topology. 
Moreover, $S$ is said to be unital if it contains the identity of $\CL(H)$
\end{defn}

\medskip

Now, we can present our generalization of Theorem \ref{thm:orig-BM} (see the paragraph following Theorem \ref{thm:orig-BM}). 

\medskip

\begin{thm}\label{thm:ext-BM}
Let $S$ be an operator system.
Then $S$ is a dual operator system if and only if it is a dual operator space with a predual Banach space $V$ such that the involution on $S$ is weak$^*$-continuous and the matrix cone on $S$ is weak$^*$-closed. 
\end{thm}
\begin{proof}
It is clear that if $S$ is a dual operator system, then the required condition holds. 
Conversely, suppose that such a predual $V$ exists. 
Let $j_S$ be the evaluation map from $S$ to the $C^*$-algebra 
$$A(S):= {\bigoplus}_{n\in \BN} C(B_{M_n(S^*)}^+;M_n);$$
where $B_{M_n(S^*)}^+$ is equipped with the weak$^*$-topology. 
Since the compact Hausdorff space $\mathcal{Q}^S_n$ of all completely positive complete contractions from $S$ to $M_n$ coincides with $B_{M_n(S^*)}^+$, we know from Lemma 2.4(d) and Theorem 2.6 of \cite{Ng2} that $j_S$ is a complete isometry (as $S$ is an operator system). 

On the other hand, under the hypothesis, $S$ is a dual MOS in the sense of \cite{Ng1} with the fixed predual $D_\# = V$. 
We regard $V$ as a sub-SMOS of the SMOS $S^*$. 
Then the set $\mathcal{WQ}^S_n$ of all $\sigma(S,V)$-continuous completely positive complete contractions from $S$ to $M_n$ will coincide with $B_{M_n(V)}^+$. 
Let $\mu_{S}$ be the evaluation map from  $S$ to the von Neumann algebra 
$$N(S):= {\bigoplus}_{n\in \BN} \ell^\infty(B_{M_n(V)}^+;M_n).$$
By \cite[Proposition 3.6(b)]{Ng1}, it suffices to show that $\mu_S$ is a complete isometry. 

Indeed, we know from \cite[Corollary 3.8]{Ng1} that the canonical map from $S$ to $V^*$ is a completely isometric SMOS isomorphism. 
Hence, $S^*\cong V^{**}$ as SMOS in the canonical way. 
Consider $m\in \BN$ and $x\in M_m(S)$.
We know that 
$$\|\mu_S^{(m)}(x)\| = \sup \{\|\omega^{(m)}(x)\|: \omega\in B_{M_n(V)}^+;n\in \BN \}.$$
Thus, we know from Proposition \ref{main} that 
$$\|\mu_S^{(m)}(x)\| = \sup \{\|\varphi^{(m)}(x)\|: \varphi\in B_{M_n(V^{**})}^+;n\in \BN \} = \|j_S^{(m)}(x)\| = \|x\|.$$
This completes the proof. 
\end{proof}

\medskip

The following example tells us that the weak$^*$-closedness assumption of the matrix cone  in the above is indispensable. 

\medskip

\begin{eg}\label{eg:counter}
Let us regard the space $c_0$ of null sequences as a $C^*$-subalgebra of $\CL(\ell^2)$ in the canonical way. 
Consider $S := \CL(\ell^2)$, equipped with the matrix cone $(M_n(c_0)_+)_{n\in \BN}$. 
Clearly, $S$ is a dual operator space. 
By Lemma 2.4(d) and Theorem 2.6 of \cite{Ng2}, $S$ will be an operator system if the map $j_S$ as in the proof of Theorem \ref{thm:ext-BM} is a complete isometry; in other words, if
$$\|x\| = \sup \big\{\|\varphi^{(m)}(x)\|: \varphi\in B_{M_n(S^*)}^+;n\in \BN \big\} \qquad (x\in M_m(S); m\in \BN).$$
Now, consider $T:=\CL(\ell^2)$, equipped with the usual matrix cone $(M_n(\CL(\ell^2))_+)_{n\in\BN}$. 
Clearly, if a complete contraction $\varphi: T\to M_n$ is completely positive, then it is a completely positive map from $S$ to $M_n$. 
This shows that $B_{M_n(T^*)}^+\subseteq B_{M_n(S^*)}^+$. 
Since $T$ is an operator system,  Lemma 2.4(d) and Theorem 2.6 of \cite{Ng2} tells us that 
\begin{align*}
\|x\| = \sup \big\{\|\varphi^{(m)}(x)\|: \varphi\in B_{M_n(T^*)}^+;n\in \BN \big\}
\leq \sup \big\{\|\varphi^{(m)}(x)\|: \varphi\in B_{M_n(S^*)}^+;n\in \BN \big\}\leq \|x\|, 
\end{align*}
for any $m\in \BN$ and $x\in M_m(S) = M_m(T)$.
This shows that $S$ is an operator system. 
However, $(c_0)_+$ is definitely not weak$^*$-closed in $\CL(\ell^2)$. 
\end{eg}

\medskip

\section*{Acknowledgement}

\medskip

The authors are supported by the National Natural Science Foundation of China (11871285).

\end{document}